\documentclass[11pt,reqno,a4paper]{amsart}

\usepackage{graphicx}
\usepackage{color}
\usepackage[USenglish]{babel}
\usepackage{latexsym}
\usepackage{amsmath}
\usepackage{amsfonts}
\usepackage{amssymb}
\usepackage{esint}
\usepackage[all]{xy}

\newcommand{\grad}{\nabla}

\renewcommand{\to}{\rightarrow}
\newcommand{\pa}{\partial}
\newcommand{\diff}{\!\setminus\!}
\newcommand{\ino}{\int_{\Omega}}

\newcommand{\dx}{\, dx}

\newcommand{\ainf}{\mbox{as\;}\;n\to+\infty}

\newcommand{\rife}[1]{(\ref{#1})}
\newcommand{\ov}[1]{\overline{#1}}

\newcommand{\sscp}{\scriptscriptstyle}
\newcommand{\dsp}{\displaystyle}

\renewcommand{\dfrac}{\displaystyle\frac}
\newcommand{\finedim}{\hspace{\fill}$\square$}
\newcommand{\intbar}{\mathop{\int\makebox(-15.5,0){\rule[6pt]{.7em}{0.3pt}}\kern-6pt}\nolimits}

\newcommand{\dt}{\delta}

\newcommand{\al}{\alpha}

\newcommand{\sg}{\sigma}

\newcommand{\om}{\Omega}
\newcommand{\ww}{\widehat{w}}

\newcommand{\graf}[1]{\left\{\begin{array}{ll}#1\end{array}\right.}

\newcommand{\D }{\Delta }

\newcommand{\e}[1]{{\,\dsp e^{\dsp #1}}}

\renewcommand{\l }{\tau }

\newcommand{\n }{\nabla }

\newcommand{\s }{\sigma }

\renewcommand{\o }{\omega }
\renewcommand{\O }{\Omega }

\def\o{\omega}
\def\p{\partial}

\newcommand{\be}{\begin{equation}}
\newcommand{\ee}{\end{equation}}
\newcommand{\beq}{\begin{equation}}
\newcommand{\eeq}{\end{equation}}

\newcommand{\R}{\mathbb{R}}

\newcommand{\ci}{\mathbb{C}}

\newtheorem{theorem}{Theorem}[section]
\newtheorem{proposition}[theorem]{Proposition}
\newtheorem{definition}{Definition}[section]
\newtheorem{corollary}[theorem]{Corollary}
\newtheorem{remark}[theorem]{Remark}
\newtheorem{example}[theorem]{Example}

\newtheorem{lemma}[theorem]{Lemma}

\newcommand{\bpr}{\begin{proposition}}
\newcommand{\epr}{\end{proposition}}
\newcommand{\bex}{\begin{example}\rm}
\newcommand{\eex}{\end{example}}
\newcommand{\brm}{\begin{remark}\rm}
\newcommand{\erm}{\end{remark}}
\newcommand{\bdf}{\begin{definition}\rm}
\newcommand{\edf}{\end{definition}}
\newcommand{\bte}{\begin{theorem}}
\newcommand{\ete}{\end{theorem}}
\newcommand{\ble}{\begin{lemma}}
\newcommand{\ele}{\end{lemma}}
\newcommand{\bco}{\begin{corollary}}
\newcommand{\eco}{\end{corollary}}

\begin{document}

\newtheorem{lem}{Lemma}[section]
\newtheorem{pro}[lem]{Proposition}
\newtheorem{thm}[lem]{Theorem}
\newtheorem{rem}[lem]{Remark}
\newtheorem{cor}[lem]{Corollary}
\newtheorem{df}[lem]{Definition}

\title[On the first eigenvalue of Liouville-type problems]
{On the first eigenvalue of Liouville-type problems}

\author{
{Daniele Bartolucci},
{Paolo Cosentino}, {Aleks Jevnikar} \and
{Chang-Shou Lin}}

\address{Daniele Bartolucci, Department of Mathematics, University of Rome {\it ``Tor Vergata"} \\  Via della ricerca scientifica n.1, 00133 Roma, Italy. }
\email{bartoluc@mat.uniroma2.it}

\address{Paolo Cosentino, Department of Mathematics, University of Rome {\it ``Tor Vergata"} \\  Via della ricerca scientifica n.1, 00133 Roma, Italy. }
\email{cosentino@mat.uniroma2.it}

\address{Aleks Jevnikar, Department of Mathematics, Computer Science and Physics, University of Udine, Via delle Scienze 206, 33100 Udine, Italy.}
\email{aleks.jevnikar@uniud.it}

\address{Chang-Shou Lin, Taida Institute for Mathematical Sciences and Center for Advanced Study in Theoretical Sciences, National Taiwan University, Taipei, Taiwan.}
\email{cslin@math.ntu.edu.tw}


\thanks{D.B. and P.C. are partially supported by the MIUR Excellence Department Project\\ MatMod@TOV
awarded to the Department of Mathematics, University of Rome Tor Vergata.}

\begin{abstract}
The aim of this note is to study the spectrum of a linearized Liouville-type problem, characterizing the case in which the first eigenvalue is zero. Interestingly enough, we obtain also point-wise information on the associated first eigenfunction. To this end, we refine the Alexandrov-Bol inequality suitable for our problem and characterize its equality case.
\end{abstract}

\maketitle


\

\section{Introduction}\label{sec1}

We are concerned with subsolutions of the Liouville equation, that is
\beq\label{eq0}
-\D w \leq \e{w}  \quad  \mbox{ in } \quad \om,
\eeq

where 
$w\in C^2(\O)\cap C^{0}(\ov\O)$ and we assume that $\om\subset\R^{2}$ is an open, bounded domain. We consider the eigenvalue problem associated to (\ref{eq0}), that is,
\begin{equation} \label{lin1}
\left\{ \begin{array}{ll}
-\D\phi -e^w\phi=\hat\nu e^w\phi & \mbox{in } \O, \vspace{0.2cm}\\
 \phi=0 & \mbox{on } \p \O.
\end{array}
\right.
\end{equation}

A lot of work has been done to obtain sufficient conditions to guarantee that the
first eigenvalue $\hat\nu_1$ of \eqref{lin1} satisfies $\hat\nu_1>0$, which is in turn related to non-degeneracy and uniqueness properties of the associated Liouville problem, see for example \cite{band, bjl, BLin3, CCL, suz}. In particular, it turns out that for $\int_\O e^w\,dx\leq 4\pi$ one has $\hat\nu_1\geq0$. However, at least to our knowledge, we do not have a characterization of the case $\hat\nu_1=0$. Our aim here is to fill this gap.

\

Here and in the rest of the paper a multiply connected domain is a connected but not simply connected domain and $B_\dt=\{x\in \R^2\,:\,|x|<\dt\}$, which we will sometime identify
with $\{z\in \ci\,:\,|z|<\dt\}$. Moreover, let us set
\begin{equation} \label{U}
	U_{\l}(x) = \ln\left( \dfrac{ \l}{ 1+\frac{\l^2}{8}|x|^{2}} \right)^2,\quad \l>0,
\end{equation}
which satisfies
$$
	\D U_{\l} + e^{U_{\l}}=0 \mbox{ in }\R^2,\quad \int\limits_{\R^2}e^{U_{\l}}=8\pi.
$$
Then, we have the following.

\begin{thm} \label{pro-nodal}
Let $\om\subset\R^2$ be an open, bounded domain whose boundary is the union of finitely many rectifiable Jordan curves and
$w\in C^2(\O)\cap C^{0}(\ov\O)$ be a solution of {\rm \rife{eq0}}. Let $\hat\nu_1$ be the first eigenvalue of \eqref{lin1} and assume that $\int_\O e^w\,dx\leq 4\pi$.

\medskip

$(j)$ If $\om$ is simply connected then $\hat\nu_1\ge0$ and $\hat\nu_1=0$ happens if and only if:
\begin{itemize}
 \item[] $(a)_0$ $\int_\O e^w\,dx=4\pi$;\\
$(a)_1$ the equality sign holds in {\rm \rife{eq0}} for any $x\in\om$;\\
$(a)_2$ there exists a conformal map $\Phi:\ov{B_1} \to \ov{\om}$ such that,
$$
\e{w(\Phi(z))}|\Phi^{'}(z)|^2|dz|^2=e^{U_{\sqrt8}(z)}|dz|^2,  \quad z\in \ov{B_1};
$$
$(a)_3$ the first eigenfunction $\phi$, relative to $\hat\nu_1=0$, takes the
form $\phi(x)=\varphi(\Phi^{-1}(x))$, where $\varphi(z)=\frac{1-|z|^2}{1+|z|^2}$.
\end{itemize}

\smallskip

\textcolor{white}{$(j)$} Assume that $w$ satisfies, for $c\in\R$,
$$
w=c \quad \mbox{on } \p\O,
$$
then $\hat\nu_1=0$ happens if and only if in addition to $(a)_0,(a)_1,(a)_2,(a)_3$ it holds:
\begin{itemize}
 \item[] $(a)_4$ there exists $\delta>0$ and $\theta\in\R$ such that, up to a translation, $\Phi(z)=\dt e^{i\theta} z$ and then $\om=B_{\dt}$ and
$$
\phi(x)=\varphi(\dt^{-1}x).
$$\end{itemize}

\medskip

$(jj)$ If $\om$ is a multiply connected domain of class $C^1$ and $w\in C^{1}(\ov\O)$ satisfies, for $c\in\R$
\begin{align*}
&w\ge c  \quad  \mbox{ in } \O, \\
&w=c  \quad  \mbox{ on } \p \om,
\end{align*}
then $\hat\nu_1> 0$.
\end{thm}

\medskip

\brm\label{remeig}{\it From a geometric viewpoint, in the case $\hat\nu_1=0$, $(\om,e^{w}|dx|^2)$ must be conformally equivalent to a hemisphere of the sphere of radius $\sqrt{2}$, say $\mathbb{S}_{\sscp \sqrt{2}}$. Actually, when $w=c$ on $\pa\om$ we further show that the conformal map is proportional to the identity. This in turn gives new point-wise information on the associated first eigenfunction. We refer to Remark \ref{remal} for further geometric interpretations.}
\erm

The proof of Theorem \ref{pro-nodal} is based on a refined version of the Alexandrov-Bol inequality suitable to be applied to subsolutions of the Liouville equation \eqref{eq0}, possibly on multiply connected domains, and a careful analysis of its equality case. Indeed, such isoperimetric inequality plays a major role in the symmetrization argument needed to estimate the first eigenvalue.

\

We conclude the introduction with the following remark about the singular counterpart of \eqref{eq0}.

\brm{\it We briefly address here what happens in case we consider the more general Liouville problem
$$
	-\D w \leq h(x)\e{w} \ \mbox{a.e. in } \O,
$$
where $h$ is a positive weight, possibly singular. Suppose for the moment $\O$ is simply connected. Following the arguments in \cite{bjl}, it would turn out that if $\log(h)$ is subharmonic in $\O$ then we always have $\hat\nu_1> 0$ for $\int_\O h(x)e^w\,dx\leq 4\pi$. On the other hand, if $\log(h)$ is superharmonic then we would have a modified sharp threshold $\int_\O h(x)e^w\,dx= 4\pi(1-\alpha)$, for some $\alpha>0$ depending on $h$, and it should be possible to classify the first eigenfunction, relative to $\hat\nu_1=0$, in terms of $\varphi_a(z)=\frac{1-|z|^{2(1-\alpha)}}{1+|z|^{2(1-\alpha)}}$. The case of $\O$ multiply connected is more subtle, as discussed in \cite{BLin3}, and has not been studied in full generality. We therefore postpone this discussion to a future work.}
\erm

\

The paper is organized as follows. In section \ref{sec2} we discuss the Alexandrov-Bol inequality on simply connected domains, characterizing the equality case. The case of multiply connected domains is treated in section \ref{sec3}. Finally, section \ref{sec4} is devoted to the proof of Theorem \ref{pro-nodal}. Part of the proof of the Alexandrov-Bol inequality on multiply connected domains is postponed to the appendix.

\bigskip
\bigskip

\section{The Alexandrov-Bol inequality on simply connected domains}\label{sec2}

In this section we introduce the Alexandrov-Bol inequality on simply connected domains, which was first derived in the analytical framework in \cite{band} and later generalized in \cite{suz}. We further refine the argument of \cite{suz}, giving a full characterization of the equality case.

\begin{pro}\label{bolold}
Let $\O\subset\R^2$ be a open, bounded and simply connected domain whose boundary is a rectifiable Jordan curve and
$w\in C^2(\O)\cap C^{0}(\ov\O)$ be a solution of {\rm \rife{eq0}} which satisfies,
$$
\int_{\O} \e{w}\leq 8\pi.
$$
Let $\omega\subseteq \O$ be any open subset whose boundary is the union of finitely many rectifiable Jordan curves. Then it holds:
\begin{equation} \label{bol}
	\left( \int_{\p\omega}\left(e^{w}\right)^{\frac 12}\,d\s \right)^2 \geq
\frac 12 \left( \int_\omega e^w\,dx \right)\left( 8\pi-\int_\omega e^w\,dx \right).
\end{equation}
Moreover, the equality holds in {\rm \rife{bol}} if and only if:\\
\begin{itemize}
 \item[] $(i)_1$ the equality sign holds in {\rm \rife{eq0}} for any $x\in\o$;\\
$(i)_2$ $\o$ is simply connected;\\
$(i)_3$ there exists $\l>0$ and a conformal map $\Phi:\ov{B_1} \to \ov{\o}$ such that,
\beq\label{equality1}
\e{w(\Phi(z))}|\Phi^{'}(z)|^2|dz|^2=e^{U_\l(z)}|dz|^2,  \quad z\in \ov{B_1}.
\eeq
\end{itemize}

In particular, if either the inequality in {\rm \rife{eq0}} is not an equality on $\o$ or  $\o$ is not simply connected, then the inequality in {\rm \rife{bol}} is always strict.

\smallskip

Assume that $w$ satisfies, for $c\in\R$,
$$
w=c \quad \mbox{on } \p\o,
$$
then the equality holds in {\rm \rife{bol}} if and only if, in addition to $(i)_1$-$(i)_2$-$(i)_3$, it holds,\\
\begin{itemize}
 \item[] $(i)_4$ there exists $\dt>0$ and $\theta\in\R$ such that, up to a translation, $\Phi(z)=\dt e^{i\theta}z$ and then in particular $\o=B_{\dt}$ and $w(x)=U_{\l\delta^{-1}}(x)$ .
 \end{itemize}

\end{pro}

\medskip

\brm\label{remal}{\it
There is a well known geometric meaning behind the result, see {\rm \cite{Ban}} and references therein.
In particular, in the equality case, $(i)_3$ speaks that the abstract surface $(\o,e^{w}|dx|^2)$
is conformally equivalent to $(B_1,e^{U_\l(z)}|dz|^2)$.
Actually $e^{U_{\l}(z)}|dz|^2$ is the local expression, after stereographic projection,
of the standard metric of $\mathbb{S}_{\sscp \sqrt{2}}$.
Whence, in particular the equality holds if and only if $(\o,e^{w}|dx|^2)$ is conformally equivalent to
a geodesic disk on $\mathbb{S}_{\sscp \sqrt{2}}$. Interestingly enough, our result shows that if $w=c$ on $\pa\o$,
then the equality holds if and only if $(\o,e^{w}|dx|^2)$ coincides with the local coordinates expression
of a geodesic disk, that is, the conformal map is proportional to the identity in this case.}
\erm
\bigskip

\proof[Proof of Proposition \ref{bolold}]$\,$\\
We first prove the proposition in case $\o\subseteq \O$ is simply connected, whence $\pa \o$ will be a rectifiable Jordan curve.
Let
$$
f:=-\Delta w-\e{w}\leq 0 \mbox{ in }\o,
$$
and $h_-$ be the unique solution of
$-\D h_-=f \mbox{ in }\o,\, h_-=0 \mbox{ on } \pa\o$. Next, let $h_0$ be the harmonic lifting
of $w$ on $\p \o$, that is $\D h_0=0$ in $\o$, $h_0=w$ on $\p\o$. Since $w\in C^2(\O)\cap C^{0}(\ov\O)$, then, by standard elliptic
theory (\cite{GT}), $h_-$ and $h$ are unique, $h_-$ is a subharmonic function of class $C^2(\o)\cap C^{0}(\ov\o)$
and $h_0\in C^2(\o)\cap C^{0}(\ov\o)$.
To simplify the notations let us set,
$$
h=h_0+h_-\mbox{ in }\o.
$$
At this point we define $u=w-h$, which satisfies,
\beq\label{ineq}
-\D u= e^{h}e^{u} \mbox{ in }\o,\quad u=0 \mbox{ on } \p\o.
\eeq
Clearly $u\in C^2(\o)\cap C^{0}(\ov\o)$, $u\geq 0$ in $\o$ and we define,
$$
\o(t)=\{x\in\o\,:\,u >t\}, \quad \gamma(t)=\{x\in\o\,:\,u=t\}, \quad t\in [0, t_+],
$$
where $t_+=\max\limits_{\ov \o} u$, and

$$
m(t)=\int\limits_{\o(t)}e^{h}e^{u}dx,\qquad \mu(t)=\int\limits_{\o(t)}e^{h}dx.
$$
Since $|\Delta u|$ is bounded and in particular is bounded below away from zero in $\ov{\o}$,
then it is not difficult to see that actually $m(t)$ and $\mu(t)$ are continuous in $[0, t_+]$. 
Moreover, we notice that, by well-known arguments, the set $\{x\in\o|\grad u(x)=0\}\cap u^{-1}([0,t_+])$ is of measure zero and we can use the co-area formula in (\cite{BW}) to deduce that $m(t)$ and $\mu(t)$ are absolutely continuous in $[0,t_+]$.
In particular the level sets have vanishing two dimensional area $|\gamma(t)|=0$ for any $t$, and we will use the fact that,
$$
m(0)=\int\limits_{\o}e^{h}e^{u}dx=\int\limits_{\o}e^{w}dx\leq 8\pi, \quad m(t_+)=0, \quad \mu(t_+)=0.
$$

By the co-area formula and the Sard Lemma we have,
\begin{equation}\label{diffm}
-m^{'}(t)=\int\limits_{\gamma(t)}\frac{e^{h}e^{u}}{|\nabla u |}d\sigma=
e^{t}\int\limits_{\gamma(t)}\frac{e^{h}}{|\nabla u |}d\sigma=
e^{t}(-\mu^{'}(t)),
\end{equation}
for a.a. $t\in [0, t_+]$, and, in view of \rife{ineq},
\begin{equation}\label{diffem}
m(t)=-\int\limits_{\o(t)}\D u=\int\limits_{\gamma(t)}|\nabla u|,
\end{equation}
for a.a. $t\in [0, t_+]$. By the Schwarz inequality we find that,
$$
-m^{'}(t)m(t)=
\int\limits_{\gamma(t)}\frac{e^{h}e^{u}}{|\nabla u |}d\sigma\int\limits_{\gamma(t)}|\nabla u|d\sigma=
e^{t}\int\limits_{\gamma(t)}\frac{e^{h}}{|\nabla u |}d\sigma\int\limits_{\gamma(t)}|\nabla u|d\sigma\geq
$$
\beq\label{hub2}
e^{t}\left(\int\limits_{\gamma(t)}{e^{h/2}} d\sigma\right)^2\geq 4\pi e^{t}\mu(t),
\mbox{ for a.a. } t\in [0, t_+],
\eeq
where in the last inequality, since $h$ is subharmonic, we used a generalization of a classical isoperimetric inequality due to Huber (\cite{hu}), which was proved in the case of open and simply connected domains. If $\o(t)$ is multiply connected the inequality is strict. For simplicity, let assume $\o(t)=\o_1(t)\diff\ov{\o_2(t)}$, with $\o_2(t)\subset\o_1(t)\subset\o$ open and simply connected domains and also $\p\o(t)=\p\o_1(t)\cup\p\o_2(t)$. We notice that $h$ is well-defined and subharmonic in both $\o_1(t)$ and $\o_2(t)$. Then we can apply Huber's inequality to both domains:
$$
\left(\int\limits_{\p\o(t)}{e^{h/2}} d\sigma\right)^2 > \left(\int\limits_{\p\o_1(t)}{e^{h/2}} d\sigma\right)^2 +\left(\int\limits_{\p\o_2(t)}{e^{h/2}} d\sigma\right)^2\geq
$$
$$
\geq 4\pi\left(\int\limits_{\o_1(t)}{e^{h}} \dx + \int\limits_{\o_2(t)}{e^{h}} \dx\right)>\int\limits_{\o(t)}{e^{h}} \dx \dx
$$
The general case follows by induction on the number of "holes". Finally, considering the case in which $\o(t)$ is not connected, we can assume for simplicity $\o(t)=\o_1(t)\cup\o_2(t)$ with $\o_1(t),\o_2(t)\subseteq\o$ open and connected subsets, $\o_1(t)\cap\o_2(t)=\varnothing$ and $\p\o(t)=\p\o_1(t)\cup\p\o_2(t)$, and make the same calculations as done before.\\
Therefore we conclude that,
\begin{equation}\label{diff1}
\frac{1}{8\pi}(m^2(t))^{'}+e^t\mu(t)\leq 0,\mbox{ for a.a. } t\in [0, t_+].
\end{equation}

In particular, because of \eqref{diffm}, we conclude that,
$$
\left(\frac{1}{8\pi}m^2(t)-m(t)+e^t\mu(t)\right)^{'}=\frac{1}{8\pi}(m^2(t))^{'}+e^t\mu(t)\leq 0,\mbox{ for a.a. } t\in [0, t_+].
$$

However, as mentioned above, the quantity in the parentheses in the l.h.s. of this inequality is continuous and absolute continuous in $[0, t_+]$, and then we also conclude that,
$$
\frac{1}{8\pi}m^2(0)-m(0)+\mu(0)\geq \frac{1}{8\pi}m^2(t_+)-m(t_+)+e^{t_+}\mu(t_+)=0,
$$
that is, by using once more the Huber (\cite{hu}) inequality,
\beq\label{hub1}
\left(\,\int\limits_{\p \o}\left({e^{w}}\right)^{\frac12} d\sigma\right)^2=
\left(\,\int\limits_{\p \o}{e^{h/2}} d\sigma\right)^2\equiv
\left(\,\int\limits_{\gamma(0)}{e^{h/2}} d\sigma\right)^2\geq 4\pi\left(\,\int\limits_{\o(0)}{e^{h}} dx\right)=
\eeq
$$
4\pi \mu(0)\geq \frac{1}{2}(8\pi m(0)-m^2(0))=\frac{1}{2}\left( 8\pi - \int\limits_{\o}e^{w}dx\right)\int\limits_{\o}e^{w}dx,
$$
which is \rife{bol}. Therefore, \rife{bol} holds for $\o$ simply connected and we now characterize the equality sign.
We first recall that the equality holds in the Huber (\cite{hu}) inequality used in
\rife{hub1} if and only if there exists $\al\in\R$ such that
$e^{h(\eta)}=e^{-\al}|\Psi^{'}(\eta)|^2$, $\eta\in \o$, where $\Psi:\o\to B_{1}$ is univalent and conformal. Since $\p\o$
is simple, then in particular by the Carath\'eodory Theorem (\cite{pom}) $\Psi$ is continuous on $\ov{\o}$  and maps one to one $\pa \o$ onto $\pa B_{1}$. Let
$\Phi=\Psi^{-1}:\ov{B_{1}}\to\ov{\o}$ and set,
$$
\xi(z)=u(\Phi(z))-\al,
$$
then we see that $\xi(z)$ satisfies,
$$
-\D \xi=e^{\xi}\mbox{ in }B_1,\quad \xi=-\al \mbox{ on }\pa B_1,
$$
and therefore it is radial (\cite{gnn}) and it is well known (\cite{suz}) that it
takes the form $\xi(z)=\xi(|z|)=U_{\l}(z)$ for some $\l>0$ which solves $U_{\l}(1)=-\al$.
Also, we deduce that,
\beq\label{sharpc}
\e{w(\Phi(z))}=\e{h(\Phi(z))}\e{u(\Phi(z))}=|\Phi^{'}(z)|^{-2}e^{U_{\l}(z)},\quad z\in \ov{B_{1}},
\eeq
which is \rife{equality1}.
Since $\log(|\Phi^{'}(z)|^2)$ is harmonic we also have,
$$
-\D_{z} w(\Phi(z))=-\D_{z} U_{\l}(z)=e^{U_{\l}(z)}=|\Phi^{'}(z)|^{2}\e{w(\Phi(z))},\quad z\in B_{1},
$$
which implies,
$$
-\D w = \e{w}\quad x\in \o.
$$
In view of the definition of $f$ and since $w\in C^{2}(\om)$, we deduce that necessarily $f\equiv 0$ and in particular,
that necessarily the equality sign in \rife{eq0} holds in $\o$.
\brm\label{remsharp}
{\it Since $\Phi$ maps $\pa B_1$ one to one and onto $\pa \o$, if $w=c$, $c\in \R$ on $\pa\o$, then $w(\Phi(z))=c$ on
$\pa B_1$. As a consequence we see from {\rm \rife{sharpc}} that $|\Phi^{'}(z)|^{2}$ is constant on $\pa B_1$ and in particular
we conclude that $\log(|\Phi^{'}(z)|^{2})$ is harmonic in $B_{1}$ and constant on $\pa B_{1}$. Therefore, we have, up to a translation,
$\Phi(z)=\dt e^{i\theta} z$, for some $\dt>0$ and $\theta\in\R$, $\o=B_{\dt}$
and in particular we infer from {\rm \rife{sharpc}} that,
$$
w(x)=U_{\l}(\dt^{-1}x)-2\log(\dt)=U_{\l \dt^{-1}}(x).
$$
It follows that in this case we have $\o=B_{\dt}$  and $w(x)=U_{\l\delta^{-1}}(x)$ for some $\dt>0$ and $\l>0$.
}

\erm

Therefore, we have shown that if $\o$ is simply connected and
the equality sign holds in \rife{bol}, then necessarily the equality sign holds in \rife{eq0} in $\o$
and \rife{equality1} holds as well. On the other side if these conditions are satisfied we have,
$$
\int\limits_{\p\omega}\left(\e{w}\right)^{\frac 12}\,d\s=
\int\limits_{\p B_\dt}\left(|\Phi^{'}(z)|^{-2}e^{U_{\l}(z)}\right)^{\frac 12}|\Phi^{'}(z)|\,d\sg(z)=
\int\limits_{\p B_\dt}e^{U_{\l}(z)/2}d\sg(z)=
\dfrac{ 2\pi\dt\l}{ 1+\frac{\l^2}{8}|\dt|^{2}},
$$
and similarly,
$$
\int_\o \e{w}=\int\limits_{ B_\dt}e^{U_{\l}}=\int\limits_{ B_\dt}\left(\dfrac{ \l}{ 1+\frac{\l^2}{8}|z|^{2}}\right)^2=
\dfrac{ \pi\dt^2\l^2}{ 1+\frac{\l^2}{8}\dt^{2}},
$$
and so we readily conclude that,
$$
\left( \int_{\p\omega}\left(e^{w}\right)^{\frac 12}\,d\s \right)^2 =
\left(\, \int\limits_{\p B_\dt}e^{U_{\l}(z)/2}d\sg(z) \right)^2 =
\frac 12 \left( \int\limits_{B_\dt}e^{U_{\l}} \right)\left( 8\pi-\int\limits_{B_\dt}e^{U_{\l}} \right)=
$$
$$
\frac 12 \left( \int_\omega e^w\,dx \right)\left( 8\pi-\int_\omega e^w\,dx \right).
$$
Therefore those conditions are necessary and sufficient as far as $\om$ is simply connected and
$\o$ is simply connected. If $\om$ is simply connected and $\o$ is an open multiply connected subdomain whose boundary is
the union of finitely many Jordan curves, as observed in \cite{CCL},
one can use the assumption $\ino e^{w}\leq 8\pi$ and work out an
induction on the number of "holes" of $\o$ , which starts by writing
$\o$ as the difference of two simply connected domains, see also the proof of Lemma \ref{lemhat} below for further details. In particular, it turns out
that the inequality \rife{bol} is always strict in this case. Finally, we consider the case of $\o$ not connected. Also in this case the proof works out by induction on the numbers of connected components and by the same calculations of Lemma \ref{lemhat} we have the strict inequality in \rife{bol}.

\finedim

\brm\label{rema2}{\it
We notice that we do not need the assumption $\int_{\O} \e{w}\leq 8\pi$ in the case $\o$ is simply connected.}
\erm

\bigskip
\bigskip

\section{The Alexandrov-Bol inequality on multiply connected domains}\label{sec3}

In this section we discuss the Alexandrov-Bol inequality on multiply connected domains, which was first derived for solutions $w$ of the Liouville equation such that $w=c$ on $\p\O$, with $c\in\R$. We refine here the argument by treating subsolutions of the Liouville equation and characterizing the equality case.

\begin{pro}\label{bolnew}
Let $\O\subset\R^2$ be a open, bounded and multiply connected domain of class $C^1$ and
$w\in C^2(\O)\cap C^{1}(\ov\O)$ be a solution of {\rm \rife{eq0}} which satisfies, for $c\in\R$,
\beq\label{eq}
w\ge c  \quad  \mbox{ in } \O,
\eeq
\beq\label{eq1}
w=c  \quad  \mbox{ on } \p \om,
\eeq
$$
\int_{\O} \e{w}\leq 8\pi.
$$
Let $\o\Subset\om$ be any relatively compact open subset whose boundary is the union of finitely many rectifiable Jordan curves. Then it holds:
\begin{equation} \label{bol2}
	\left( \int_{\p\omega}\left(e^{w}\right)^{\frac 12}\,d\s \right)^2 \geq
\frac 12 \left( \int_\omega e^w\,dx \right)\left( 8\pi-\int_\omega e^w\,dx \right).
\end{equation}
Moreover, the equality holds in {\rm \rife{bol2}} if and only if $(i)_1$-$(i)_2$-$(i)_3$ (and $(i)_4$ in case $w=c\in\R$ on $\p\o$) of Proposition \ref{bolold} hold true.

\medskip

In particular, if either the inequality in {\rm \rife{eq0}} is not an equality on $\o$ or  $\o$ is not simply connected, then the inequality in {\rm \rife{bol2}} is always strict.
\end{pro}

\medskip

\brm\label{rema3}{\it
In the case $w$ satisfies (\ref{eq0}), (\ref{eq1}) and it is also superharmonic, we can actually recover the hypothesis (\ref{eq}) by applying the weak minimum principle.
}
\erm
\bigskip

\proof[Proof of Proposition \ref{bolnew}]$\,$\\
If $\o\Subset \om$ is relatively compact and simply connected, then we can apply Proposition {\rm \ref{bolold}}, including the characterization of the equality sign.

We will show now that if $\o\Subset \om$ is relatively compact and multiply connected, then \rife{bol2} holds with the strict inequality. \\
We will denote by $\ov{\om^*}$
the closure of the union of the bounded components of $\R^2\setminus\pa\om$ and with $\om^*=\ov{\om^*}\setminus \pa(\ov{\om^*})$.
Clearly $\om\subseteq \om^*$ and $\om\equiv \om^*$ if and only if $\om$ is simply connected. Also, there is no loss of generality
in assuming $c=0$. Indeed, if $c\neq 0$, we could define,
$$
w_c(x)=w(e^{-\frac{c}{2}}x)-c, \quad x\in e^{\frac{c}{2}}\om,
$$
which satisfies \rife{eq0},\rife{eq} in $e^{\frac{c}{2}}\om$ and \rife{eq1} on $\p\big(e^{\frac{c}{2}}\om\big)$ with $c=0$, while the integrals involved in the inequality
\rife{bol2} are invariant. Therefore we assume in the rest of this proof that,
$$
c=0,
$$
and define,
$$
\widehat{w}(x)=\graf{w(x)\quad & x\in \om\\ 0 \quad\quad\quad & x \in \om^*\setminus \om.}
$$
\ble \label{lemsol}
The function $\ww$ is a solution of
\beq\label{eq01}
-\D \ww \leq \e{\ww}  \quad  \mbox{ in } \quad \om^*,
\eeq
in the sense of distributions.
\ele
\proof
Indeed, for any $\varphi\in C^{2}_0(\ov{\om^*})$, $\varphi\geq 0$ in $\om^*$, we have,
$$
-\int\limits_{\om^*}(\Delta \varphi)\ww=-\int\limits_{\om}(\Delta \varphi)w=-\int\limits_{\om}\varphi(\Delta w)+
\int\limits_{\pa \om}\varphi(\pa_{\nu} w),
$$
where $\nu$ is the exterior unit normal. Since $w\geq 0$ in
$\om$, $w=0$ on $\pa\om$ and $\om$ is of class $C^1$, $\pa_{\nu}w$ is well defined on $\pa\om$ and $\pa_{\nu} w\leq 0$ on $\pa\om$.
Therefore, since $\varphi\geq 0$, we conclude that,
$$
-\int\limits_{\om^*}\left((\Delta \varphi)\ww+\varphi e^{\ww}\right)=
-\int\limits_{\om}\varphi(\Delta w)-\int\limits_{\om^*}\varphi e^{\ww}+\int\limits_{\pa \om}\varphi(\pa_{\nu} w)\leq
$$
$$
-\int\limits_{\om}\varphi(\Delta w+e^{w})\leq 0,
$$
as claimed.
\finedim

\bigskip
\bigskip

Since $\ww$ is only Lipschitz we cannot apply directly Proposition \ref{bolold}. However in this situation
\rife{bol2} still holds whenever $\o\subseteq \o_0\Subset\om^*$, $\o_0$ is simply connected and
$\int\limits_{\o_0}\e{\ww}\leq 8\pi$.
Let us define,
$$
\widehat{\ell}(\o)=\int\limits_{\p\omega}\left(\e{\ww}\right)^{\frac 12}d\sg,\quad \widehat{m}(\o)=\int_{\o}\e{\ww}dx,
$$
then we have,
\ble\label{lemhat} Let $\ww$ be a Lipschitz continuous solution of {\rm \rife{eq01}} in the sense of distributions and let
$\o_0\Subset\om^*$ be a simply connected and relatively compact subdomain such that
$\int\limits_{\o_0}\e{\ww}\leq 8\pi$. Let $\o\subseteq\o_0$ be any open and bounded subset whose boundary is the union
of finitely many rectifiable Jordan curves. Then,
\beq\label{bolhat}
\hat{\ell}^2(\o)\geq \frac 12 \hat{m}(\o)\left( 8\pi-\widehat{m}(\o) \right),
\eeq
holds and if $\o$ is not simply connected, then the inequality is strict.
\ele
\proof
By a standard approximation argument, see Lemma 2 in \cite{BLin3},
the fact that \rife{bolhat} holds follows from the inequality \rife{bol} for $C^2(\o_0)\cap C^{0}(\ov\o_0)$ functions.
If $\o$ is connected but not simply connected we can follow the argument
in \cite{CCL} and conclude that the inequality in \rife{bolhat} is strict. Indeed, assume for simplicity that $\o=\o_1\setminus \ov{\o_2}$, $\pa\o=\pa\o_1\cup\pa\o_2$ where $\o_1$ and
$\o_2$ are open and simply connected. Then, since $\o_1=\o\cup\ov{\o_2}$, by \rife{bolhat} we have,
$$
2\hat{\ell}^2(\o)>2(\hat{\ell}^2(\pa\o_1)+\hat{\ell}^2(\pa\o_2))\geq
\hat{m}(\o\cup\o_2)(8\pi-\hat{m}(\o\cup{\o_2}))+\hat{m}({\o_2})(8\pi-\hat{m}({\o_2}))=
$$
$$
(\hat{m}(\o)+\hat{m}(\o_2))(8\pi-\hat{m}(\o)-\hat{m}(\o_2))+\hat{m}({\o_2})(8\pi-\hat{m}({\o_2}))=
$$
$$
\hat{m}(\o)(8\pi-\hat{m}(\o))+2\hat{m}(\o_2)(8\pi-\hat{m}(\o)-\hat{m}(\o_2))\geq \hat{m}(\o)(8\pi-\hat{m}(\o)),
$$
where we used $\hat{m}(\o)+\hat{m}(\o_2)=\hat{m}(\o_1)\leq 8\pi$. \\
Finally, we consider the case in which $\o$ is not connected. For simplicity, we can assume $\o=\o_1\cup\o_2$ with $\o_1,\o_2\subseteq\o_0$ open and connected subsets, $\o_1\cap\o_2=\varnothing$ and $\p\o=\p\o_1\cup\p\o_2$. Then by \rife{bolhat} we have,
$$
2\hat{\ell}^2(\o)>2(\hat{\ell}^2(\pa\o_1)+\hat{\ell}^2(\pa\o_2))\geq
\hat{m}(\o_1)(8\pi-\hat{m}(\o_1))+\hat{m}({\o_2})(8\pi-\hat{m}({\o_2}))=
$$
$$
8\pi\hat{m}(\o_1)+8\pi\hat{m}(\o_2)-(\hat{m}({\o_1})^2+\hat{m}({\o_2})^2)=
$$
$$
8\pi\hat{m}(\o)-(\hat{m}({\o_1})+\hat{m}({\o_2}))^2+2\hat{m}(\o_1)\hat{m}(\o_2)\geq
$$
$$8\pi\hat{m}(\o)-\hat{m}({\o})^2=\hat{m}(\o)(8\pi-\hat{m}(\o)).
$$
\finedim

\bigskip
\bigskip

Now, by using Lemma \ref{lemsol} and Lemma \ref{lemhat}, the validity of \eqref{bol2} with the strict inequality can be worked out following the arguments of \cite{BLin3} with minor modifications. Since this argument is not
well known, to be self-contained and for reader's convenience we carry it out in full details in the appendix. Finally, we consider the case of $\o\Subset \om$ relatively compact and not connected. Also in this case the proof works out by induction on the numbers of connected components and by the same calculations of Lemma \ref{lemhat} we have the strict inequality in \rife{bol2}.

\bigskip
\bigskip

\section{On the first eigenvalue}\label{sec4}

In this section we prove the main result about the first eigenvalue of the Liouville-type problem \eqref{lin1}.

\proof[Proof of Theorem \ref{pro-nodal}]$\,$\\
To avoid repetitions we work out the proof of $(j)$ and $(jj)$ at once.
Clearly the first eigenvalue and eigenfunction $(\hat\nu_1,\phi)$ of \eqref{lin1} correspond to the first
eigenvalue and eigenfunction $(\nu_1=\hat\nu_1+1,\phi)$ of,
\begin{equation} \label{lin3}
\left\{ \begin{array}{ll}
-\D\phi =\nu_1 e^w\phi & \mbox{in } \O, \vspace{0.2cm}\\
 \phi=0 & \mbox{on } \p \O.
\end{array}
\right.
\end{equation}
We recall that a nodal domain for $\phi\in C^0(\ov\O\,)$ is the maximal connected
component of a subdomain where $\phi$ has a definite sign. Since $\phi$ has only one nodal domain  we assume w.l.o.g. that $\phi\geq0$ in $\O$. In particular, by the maximum principle we have $\phi>0$ in $\O$. Recalling \eqref{U} we set $U(x):=U_{1}(x)$, i.e.
\begin{equation} \label{Ua}
	U(x) = \ln\left( \dfrac{ 1}{ 1+\frac{1}{8}|x|^{2}} \right)^2,
\end{equation}
which satisfies,
$$
	\D U + e^{U}=0 \mbox{ in }\R^2.
$$
Next, let $t_+=\max_{\ov\O}\phi$ and for $t>0$ let us define $\O_t=\{x\in\O \,:\,\phi>t\}$ and $R(t)>0$ such that
$$
	\int_{B_R(t)}e^{U}\,dx=\int_{\O_t} e^w\,dx.
$$
Since $\phi>0$ in $\O$ we put $\O_0=\O$ and set $R_0=\lim_{t\to0^+}R(t)$. Clearly $\lim_{t\to (t_+)^-}R(t)=0$.
Then $\phi^*:B_{R_0}\to\R$, which for $y\in B_{R_0}$, $|y|=r$, is defined by,
$$
	\phi^*(r)=\sup\{t\in(0,t_+)\,:\,R(t)>r\},
$$
is a radial, decreasing, equimeasurable rearrangement of $\phi$ with respect to the measures $e^w\,dx$ and $e^{U}dx$, and hence, in particular,
\begin{align}
	&B_{R(t)}=\{x\in\R^2\,:\,\phi^*(x)>t\}, \nonumber\\
	&\int_{\{\phi^*>t\}}e^{U}\,dx=\int_{\O_t} e^w\,dx \quad t\in[0,t_+), \nonumber\\
	&\int_{B_{R_0}}e^{U}|\phi^*|^2\,dx=\int_{\O} e^w|\phi\,|^2\,dx. \label{rearr}
\end{align}
Well known arguments (see for example \cite{bc}) show that $\phi^*$ is continuous and locally Lipschitz.
Then, by the Sard lemma, we can apply the Cauchy-Schwartz inequality and the co-area formula to conclude that,
\begin{align}
	\int_{\{\phi=t\}} |\n \phi|\,d\s & \geq \left( \int_{\{\phi=t\}} \left(e^{w}\right)^{\frac 12}\,d\s \right)^2
\left( \int_{\{\phi=t\}} \dfrac{e^w}{|\n \phi|}\,d\s \right)^{-1} \label{c-s}\\
	  & = \left( \int_{\{\phi=t\}} \left(e^{w}\right)^{\frac 12}\,d\s \right)^2 \left( -\dfrac{d}{dt}\int_{\O_t}e^w\,dx \right)^{-1}, \nonumber
\end{align}
for a.e. $t$. Next, under the assumptions either of part $(j)$ or of part $(jj)$, we are allowed to apply
the Alexandrov-Bol inequality \rife{bol},
\begin{align*}
	&\left( \int_{\{\phi=t\}} \left(e^{w}\right)^{\frac 12}\,d\s \right)^2 \left( -\dfrac{d}{dt}\int_{\O_t} e^w\,dx \right)^{-1} \\
	&\geq \frac 12 \left( \int_{\O_t} e^w\,dx \right)\left( 8\pi-\int_{\O_t} e^w\,dx \right)
\left( -\dfrac{d}{dt}\int_{\O_t} e^w\,dx \right)^{-1},
\end{align*}
for a.e. $t$. Since $\phi^*$ is an equimeasurable rearrangement of $\phi$ with respect to the measures $e^w \,dx$, $e^{U}\,dx$, and since
$e^{U}$ realizes the equality in \rife{bol}, we also conclude that,
\begin{align*}
&\frac 12 \left( \int_{\O_t} e^w\,dx \right)\left( 8\pi-\int_{\O_t} e^w\,dx \right) \left( -\dfrac{d}{dt}\int_{\O_t} e^w\,dx \right)^{-1} \\
&=\frac 12 \left( \int_{\{\phi^*>t\}} e^{U}\,dx \right)\left( 8\pi-\int_{\{\phi^*>t\}} e^{U}\,dx \right) \left( -\dfrac{d}{dt}\int_{\{\phi^*>t\}} e^{U}\,dx \right)^{-1} \\
& =\left( \int_{\{\phi^*=t\}} \left(e^{U}\right)^{\frac 12}\,d\s \right)^2 \left( -\dfrac{d}{dt}\int_{\{\phi^*>t\}} e^{U}\,dx \right)^{-1} \\
& = \int_{\{\phi^*=t\}} |\n \phi^*|\,d\s,
\end{align*}
where in the last equality we used once more the co-area formula. Therefore, we have proved that,
$$
\int_{\{\phi^*=t\}} |\n \phi^*|\,d\s \leq\int_{\{\phi=t\}} |\n \phi|\,d\s,
$$
for a.e. $t$, which in turn yields,
\begin{equation} \label{grad}
	\int_{B_{R_0}} |\n \phi^*|^2\,dx \leq \int_\O |\n \phi|^2\,dx.
\end{equation}
In deriving \eqref{grad} we have used the co-area formula together with the fact that $\int_{B_{R(t)}} |\n \phi^*|^2\,dx$ and
$\int_{\O_t} |\n \phi|^2\,dx$ are Lipschitz continuous.
By using \eqref{rearr}, \eqref{grad} and the variational characterization of $\phi$ we deduce that,
\beq\label{c0}
	\int_{B_{R_0}} |\n \phi^*|^2\,dx-\int_{B_{R_0}}e^{U}|\phi^*|^2\,dx \leq \int_\O |\n \phi|^2\,dx-\int_{\O} e^w|\phi|^2\,dx=(\nu_1-1)\int_{\O} e^w|\phi|^2\,dx,
\eeq
Moreover, $\phi^*(R_0)=0$. Now we argue by contradiction and suppose that
$\hat\nu_1<0$, so we have that $\phi$ satisfies \rife{lin3} with $\nu_1< 1$ and we have that,
$$	
 (\nu_1-1)\int_{\O} e^w|\phi|^2\,dx<0
$$
 Therefore, we conclude that the first eigenvalue of
$(-\D-e^U)(\cdot)$ on $B_{R_0}$ with Dirichlet boundary conditions is non-positive.
Consider now $\psi(x)=\frac{8-|x|^{2}}{8+|x|^{2}}$ which satisfies,
$$
	-\D\psi-e^{U}\psi=0 \mbox{ in } \R^2, \quad \psi\in C^2_0(B_{\sqrt{8}}(0)).
$$
Since the first eigenvalue is non-positive one can deduce that $R_0>\sqrt{8}$. Moreover,
\begin{equation} \label{R0}
	8\pi\frac{R_0^{2}}{8+R_0^{2}}=\int_{B_{R_0}}e^{U}\,dx=\int_\O e^w\,dx\leq 4\pi,
\end{equation}
by assumption and hence $R_0\leq\sqrt{8}$, which yields a contradiction. So, $\nu_1\geq1$, that is $\hat\nu_1\geq0$. \\

Now suppose $\hat\nu_1=0$, that is $\nu_1=1$. From (\ref{c0}), we deduce that
$$\int_{B_{R_0}} |\n \phi^*|^2\,dx-\int_{B_{R_0}}e^{U}|\phi^*|^2\,dx \leq0.$$
Therefore, we conclude that the first eigenvalue of
$(-\D-e^U)(\cdot)$ on $B_{R_0}$ with Dirichlet boundary conditions is non-positive and, arguing as before, we conclude that $R_0\geq\sqrt{8}$. On the other hand, as pointed out in (\ref{R0}), $R_0\leq\sqrt{8}$. Hence we deduce that $R_0=\sqrt{8}$ and, in particular,
$$
\int_\O e^w\,dx=\int_{B_{\sqrt{8}}}e^{U}\,dx= 4\pi
$$
Moreover, since $R_0=\sqrt{8}$, then the first eigenvalue of $(-\D-e^U)(\cdot)$ on $B_{R_0}$ with Dirichlet boundary
conditions is $0$ and $\psi$ is its eigenfunction. In particular, from (\ref{c0}), we derive that
$$
	\int_{B_{\sqrt{8}}} |\n \phi^*|^2\,dx=\int_{B_{\sqrt{8}}}e^{U}|\phi^*|^2\,dx = \int_{\O}e^w|\phi|^2\,dx=\int_{\O} |\n \phi|^2\,dx
$$
and that all the inequalities
used to obtain \rife{c0} must be equalities. In particular, for a.e. $t\in [0,t_+)$ we have:
\begin{equation} \label{equal}
\left( \int_{\{\phi=t\}} \left(e^{w}\right)^{\frac 12}\,d\s \right)^2 =\frac 12 \left( \int_{\O_t} e^w\,dx \right)\left( 8\pi-\int_{\O_t} e^w\,dx \right),
\end{equation}
Then we can choose a sequence $t_n\to 0^+$, $\ainf$,
such that the equality (\ref{equal}) holds for any $n$. In both situations $(j)$ and $(jj)$, we can apply Proposition \ref{bolold} or Proposition \ref{bolnew} and use the characterization of the equality sign in the Alexandrov-Bol inequality. In particular, we derive that $\O_{t_n}$ are simply connected. This, together with the fact that $\O_{t_n}\Delta\O\rightarrow \varnothing$, as $n\rightarrow +\infty$, implies that $\O$ is simply connected and proves (jj).  Now, taking the $\liminf$ in (\ref{equal}) and using for example Theorem 2.3 in \cite{D}, we deduce that
$$
\left( \int_{\p \O} \left(e^{w}\right)^{\frac 12}\,d\s \right)^2 \leq \frac 12 \left( \int_{\O} e^w\,dx \right)\left( 8\pi-\int_{\O} e^w\,dx \right),
$$


and the latter inequality turns out to be an equality by \eqref{bol} in Proposition \ref{bolold} with $\o=\O$. Then, applying again the equality case of Proposition \ref{bolold}, we see from $(i)_1$ and $(i)_3$ that the equality holds in \rife{eq0}
in $\om$ and that,
$$
\e{w(\Phi(z))}|\Phi^{'}(z)|^2|dz|^2=e^{U_\l(z)}|dz|^2,\quad z\in \ov{B_{1}},
$$
where $\Phi:B_1\to \om$ is conformal and univalent. Therefore we have $\int\limits_\O e^w=\int\limits_{B_1}e^{U_{\l}}$ and in
particular,
$$
4\pi=\int\limits_\O e^w =\int\limits_{B_{R_0}}e^{U}=\int\limits_{B_1}e^{U_{\l}}=\int\limits_{B_\l}e^{U_{1}},
$$
which immediately implies that $\l=\sqrt{8}$. Therefore, we have that $\varphi(z)=\phi(\Phi(z))$ satisfies
\beq \label{lin3.1}
\left\{ \begin{array}{ll}
-\D\varphi =\nu_1 |\Phi^{'}(z)|^2 \e{w(\Phi(z))}\varphi =e^{U_{\sscp \sqrt{8}}(z)}\varphi & \mbox{ in } B_1, \vspace{0.2cm}\\
 \varphi=0 & \mbox{on } \p B_1,
\end{array}
\right.
\eeq
The function $\varphi(z)=\frac{1-|z|^2}{1+|z|^2}=\psi(\sqrt{8}z)$ is a positive solution of \rife{lin3.1} and therefore
it is the first eigenfunction for $(-\D - e^{U_{\sqrt{8}}(z)})_{B_1}$ with Dirichlet boundary conditions. We conclude
that $\phi(x)=\varphi(\Phi^{-1}(x))$ which is $(a)_3$.\\
If now we assume that $(a)_0-(a)_3$ are true, it can be easily proved that $\nu_1=1$, that means $\hat\nu_1=0$. \\
We conclude the $(j)$ part, assuming that $w=c$ on $\pa\om$. Then we can apply what specified in the Remark \ref{remsharp} which implies that, up to a translation,
$\Phi(z)=\dt e^{i\theta} z$ for some $\dt>0$ and $\theta\in\R$, and we have $\O=B_{\dt}$ and $\phi(x)=\varphi(\dt^{-1}x)$, that is $(a)_4$. Moreover, this condition, together with $(a)_0-(a)_3$, is sufficient to have $\hat\nu_1=0$. \\

\finedim

\bigskip
\bigskip

\
\section{Appendix}\label{sec5}

In this section we complete the proof of the Alexandrov-Bol inequality on multiply connected domains, i.e. Proposition \ref{bolnew}, following the strategy in \cite{BLin3}.

\medskip

We recall that $\o$ is assumed to be relatively compact and multiply connected and we aim to show \eqref{bol2} holds with the strict inequality. We start with the case where each bounded component of $\R^2\backslash \ov\o$
contains at least one bounded component of $\R^2\backslash \p \O$.\\

Let $\O_0$ be the
union of the bounded components of $\R^2\backslash \p\O$ bounded by $\p\o$. Moreover, let $\o_0$ be the union of all bounded
simply connected components of $\R^2\backslash \p \o$. Thus
$\O_0\subset \o_0$ and we let
$$
\o^*=\o_0\backslash \O_0.
$$
We have $\o^*\subset \O$. Moreover, $\o^* \cup \O_0$ is a union of simply connected domains and
$\o\cup\ov{\o^\ast\cup\O_0}$ is a simply connected domain. Denote by $\p_0\o$ the boundary of $\o^* \cup \O_0$
and $\p_1\o=\p\o\setminus\p\o^*$. Then $\p_1\o=\p(\o\cup\ov{\o^\ast\cup\O_0})$ and it holds
$$
\p \o = \p_1 \o \cup \p_0 \o.
$$

To simplify the notations we define
$$
{\ell}(\o)=\int\limits_{\p\omega}\left(\e{w}\right)^{\frac 12}d\sg,\quad {m}(\o)=\int_{\o}\e{w}dx,
$$
whenever $\o\Subset \om$.\\

\noindent {\bf Case 1.}\ \ \ $\widehat{m}(\o^* \cup \O_0) \geq 8\pi.$\\
Since $w\ge0$ in $\om$, by the isoperimetric inequality we have,
\beq\label{eqn2.8}
2{\ell}(\p_0\o)^2=2\left(\,\int\limits_{\p_0\omega}\left(\e{w}\right)^{\frac 12}d\sg\right)^2\ge2\left(\,\int\limits_{\p_0\omega}d\sg\right)^2\geq
8\pi \int\limits_{\omega^*\cup \om_0}dx>8\pi\widehat{m}(\om_0).
\eeq
\bigskip

Since $\widehat{m}(\o^*\cup\O_0)\geq 8\pi$, we have
\[
\widehat{m}(\O_0)\geq8\pi-\widehat{m}(\o^*)\equiv8\pi-m(\o^*),
\]
and then by using \eqref{eqn2.8},
\begin{equation}\label{eqn2.9}
2\ell^2(\p_0\o)>8\pi(8\pi-m(\o^*)).
\end{equation}

Similarly, since $\widehat{m}(\o\cup\ov{\o^\ast\cup\O_0}) > \widehat{m}(\o^* \cup \O_0)\geq 8\pi$,
we deduce
\beq{\label{eqn2.10}}
2\ell^2(\p_1\o) > 8\pi (8\pi-m(\o\cup \o^*)).
\eeq

Recalling \eqref{eqn2.9}, (\ref{eqn2.10}) and the fact that since $\o^*\subset \om$, then $m(\o^*)\leq 8\pi$, we have
\begin{eqnarray*}
2\ell^2(\p\o) & = & 2(l(\p_1\o)+l(\p_0\o))^2 > 2 [l^2(\p_1\o)+l^2(\p_0\o)]\\
& > & m(\o^*) (8\pi-m(\o^*)) + m(\o\cup\o^*) (8\pi -m(\o\cup\o^*))\\
& = & m(\o^*) (8\pi-m(\o^*)) + (m(\o)+m(\o^*)) (8\pi - m(\o) -m(\o^*))\\
& = & m(\o) (8\pi-m(\o)) + m(\o^*) (16\pi -2m(\o) -2m(\o^*)).
\end{eqnarray*}

Finally, since
$$
m(\o)+m(\o^*)\leq m(\O) \leq 8\pi,
$$
then we conclude that,
$$
2\ell^2(\p\o)
> m(\o) (8\pi-m(\o)),
$$
which proves that \rife{bol2} holds with the strict inequality.

\bigskip

\noindent {\bf Case 2.}\ \ \ $\widehat{m}(\o^*\cup \O_0)<8\pi$ and $\widehat{m}(\o\cup\ov{\o^\ast\cup\O_0})\geq 8\pi$.\\

Since $\widehat{m}(\o^*\cup \O_0)<8\pi$ and $\o^*\cup \O_0\Subset \om^*$ is union of simply connected domains, \rife{bolhat} yields,
\begin{eqnarray}{\label{eqn2.11}}
2\ell^2(\p_0\o) & \geq & \widehat{m}(\o^*\cup \O_0)(8\pi-\widehat{m}(\o^*\cup \O_0))\\
 & = & (m(\o^*)+\widehat{m}(\O_0)) (8\pi-m(\o^*)-\widehat{m}(\O_0)).\nonumber
\end{eqnarray}
Observing that (\ref{eqn2.8}) holds, we infer that
\beq{\label{eqn2.12}}
\ell(\p_0\o) \geq \sqrt{4\pi \widehat{m}(\O_0)}, \ \ \ \ \ell(\p_1\o) \geq \sqrt{4\pi \widehat{m}(\O_0)}.
\eeq

Observe that (\ref{eqn2.10}) still holds in this case and therefore, by (\ref{eqn2.11}),
(\ref{eqn2.12}) we have
\begin{eqnarray*}
2\ell^2(\p\o) & = & 2[\ell^2 (\p_1\o)+2\ell(\p_1\o)\ell(\p_0\o)+\ell^2(\p_0\o)]\\
& \geq & 8\pi (8\pi -m(\o)-m(\o^*)) \\
& & + (m(\o^*)+\widehat{m}(\O_0)) (8\pi-m(\o^*)-\widehat{m}(\O_0))+16\pi \widehat{m}(\O_0)\\
& = & 8\pi (8\pi -m(\o)) - m^2(\o^*) + \widehat{m}(\O_0) (24\pi -2m(\o^*)-\widehat{m}(\O_0))\\
& = & m(\o) (8\pi-m(\o))+[(8\pi-m(\o))^2 - m^2(\o^*)]\\
& &  + \widehat{m} (\O_0) (24\pi-2m(\o^*)-\widehat{m}(\O_0)).
\end{eqnarray*}

\noindent Now, since $m(\o)+m(\o^*)\leq m(\O)\leq8\pi$,
\[
m(\o^*)\leq8\pi-m(\o).
\]
On the other hand, $\widehat{m}(\o^*\cup\O_0)<8\pi$, hence
\[
2m(\o^*)+2\widehat{m}(\O_0)<16\pi.
\]

We have proved that $2\ell^2(\p\o)>m(\o)(8\pi-m(\o))$, that is \rife{bol2} holds with the strict inequality.\\

\noindent {\bf Case 3.}\ \ \ $\widehat{m}(\o\cup\ov{\o^\ast\cup\O_0})<8\pi$.

Since $\o\Subset \om$ by assumption, then $\o\cup\ov{\o^\ast\cup\O_0}\Subset \om^*$. Moreover, $\o\cup\ov{\o^\ast\cup\O_0}$ is simply connected. Therefore, Lemma \ref{lemhat} yields,
\[
2\ell(\p\o)^2\geq m(\o)(8\pi-m(\o)),
\]
that is, \rife{bol2} holds and in particular, since $\o$ is not simply connected by assumption, it holds with
the strict inequality sign.\\

To finish the proof, we are left with the case $\p\o$ bounds some simply connected subdomains $\o_1,\ldots,\o_k$ of $\O$. Clearly, $\o\cup\ov{\o_1}\cup\dots\ov{\o_k}$ is simply connected in
$\O$.
Therefore, using in this case the standard inequality \eqref{bol} we infer
$$
2\ell^2(\p(\o\cup \ov\o_1 \cup \ldots \cup \ov\o_k))\geq m(\o\cup \ov\o_1 \cup \ldots \cup \ov\o_k)[8\pi - m(\o\cup \ov\o_1 \cup \ldots \cup \ov\o_k)],
$$
\begin{equation*}
2\ell^2 (\p \o_j) \geq m(\o_j) (8\pi - m(\o_j)),\;j\in \{1,\cdots,k\}.
\end{equation*}

For $k=1$ we readily have
\begin{eqnarray*}
2\ell(\p\o)^2 & > & 2\ell (\p(\o\cup\ov\o_1))^2+2\ell(\p\o_1)^2\\
& \geq & m(\o\cup\ov\o_1) (8\pi -m(\o\cup\ov\o_1))+m(\o_1)(8\pi-m(\o_1))\\
& = & m(\o) (8\pi-m(\o))+m(\o_1) (16\pi-2m(\o_1)-2m(\o))\\
& \geq & m(\o) (8\pi-m(\o)),
\end{eqnarray*}
and we deduce that \rife{bol2} holds with the strict inequality in this case as well.
The case $k>1$ is similar.\\

\finedim



\


\begin{thebibliography}{99}

\bibitem{band}
C. Bandle, \emph{On a differential Inequality and its applications to Geometry}, Math. Zeit. \textbf{147} (1976), 253-261.

\bibitem{Ban} C. Bandle, {\sl Isoperimetric inequalities and applications}, Pitmann, London, 1980.

\bibitem{bc}
D. Bartolucci, D. Castorina,  {\em On a singular Liouville-type equation and the
Alexandrov isoperimetric inequality}, Ann. Sc. Norm. Super. Pisa Cl. Sci. {\bf (5)} XIX (2019) 1-30.

\bibitem{bjl} D. Bartolucci, A. Jevnikar, C.S. Lin,
{\em Non-degeneracy and uniqueness of solutions to singular mean field equations on bounded domains},
J.D.E. {\bf 266} (2019), no. 7, 716-741.


\bibitem{BLin3} D. Bartolucci, C.S. Lin, {\em Existence and uniqueness for
Mean Field Equations on multiply connected domains at the critical parameter},
{Math. Ann.} {\bf 359} (2014), 1-44.

\bibitem{BW} J.E. Brothers, W.P. Ziemer, {\em Minimal rearrangements of Sobolev functions}, {J. Reine Angew. Math.} {\bf 384} (1988), 153-179.

\bibitem{CCL} S.Y.A. Chang, C.C. Chen, C.S. Lin, {\em Extremal functions for a mean field equation in two dimension},
in: "Lecture on Partial Differential Equations", New Stud. Adv. Math. {\bf 2} Int. Press, Somerville, MA, 2003, 61-93.

\bibitem{D} W.R. Derrick, \emph{Weighted Convergence in length}, Pac. J. Math. {\bf 43} (1972), 307-315.

\bibitem{gnn} B. Gidas, W.M. Ni, L. Nirenberg, {\it Symmetry and related
properties via the maximum principle}, Comm. Math. Phys. {\bf 68} (1979),
209-243.

\bibitem{GT} D. Gilbarg, N. Trudinger,
 "Elliptic Partial Differential Equations of Second Order", Springer-Verlag, Berlin-Heidelberg-New York (1998).

\bibitem{hu}
A. Huber, {\em Zur Isoperimetrischen Ungleichung Auf Gekr\"ummten Fl\"achen}, Acta. Math. \textbf{97} (1957), 95-101.

\bibitem{ne} Z. Nehari, {\em On the principal frequency of a membrane}, Pac. J. Math. {\bf 8} (1958), 285-293.

\bibitem{pom} Ch. Pommerenke, {\sl Boundary Behaviour of Conformal Maps}, Grandlehren der Math.
Wissenschaften, {\bf 299}, p. 300, Springer-Verlag, Berlin-Heidelberg, 1992.

\bibitem{suz} T. Suzuki, {\em Global analysis for a two-dimensional elliptic eigenvalue problem with the exponential
                nonlinearity}, Ann. Inst. H. Poincar\'e Anal. Non Lin\'eaire {\bf 9} (1992), no. 4, 367-398.
							
\end{thebibliography}
\end{document}